\documentclass{amsart}
\usepackage{epsf}
\usepackage{amscd}

\newsymbol\rtimes 226F

\newcommand{\h}{Hoch\-schild\ }
\newcommand{\hc}{Hoch\-schild co\-chain}
\newcommand{\hg}{ho\-mot\-o\-py G-}
\newcommand{\coh}{co\-ho\-mol\-o\-gy}

\newcommand{\coder}{\operatorname{Coder}}

\newcommand{\F}[1]{\underline{\mathcal M}(#1)}
\newcommand{\FF}{\underline{\mathcal M}}

\newcommand{\Hom}{\operatorname{Hom}}
\newcommand{\id}{\operatorname{id}}
\newcommand{\IN}{{\operatorname{in}}}

\newcommand{\MM}{\mathcal{M}}
\newcommand{\oE}{\overline E^1}
\newcommand{\PP}{\mathcal{P}}
\newcommand{\pr}{\operatorname{pr}}

\newcommand{\X}[1]{\underline{\mathcal{M}}_{0,#1}}

\newcommand{\del}{\partial}
\newcommand{\nc}{{\mathbb{C}}}

\newcommand{\nr}{{\mathbb{R}}}
\newcommand{\nz}{{\mathbb{Z}}}

\newtheorem{thm}{Theorem}[section]

\newtheorem{prop}[thm]{Proposition}
\newtheorem{crl}[thm]{Corollary}
\newtheorem{conj}[thm]{Conjecture}
\newtheorem{quest}[thm]{Question}

\theoremstyle{definition}

\newtheorem{df}{Definition}[section]

\theoremstyle{remark}

\newtheorem{rem}{Remark}
\newtheorem{ack}{Acknowledgments}

\newcommand{\abs}[1]{\lvert#1\rvert}

\begin{document}

\title{Homotopy Gerstenhaber algebras}

\author{Alexander A. Voronov}
\thanks{Published in \emph{Conf\'erence Mosh\'e Flato 1999}
(G. Dito and D. Sternheimer, eds.), vol.\ 2. Kluwer Academic
Publishers, the Netherlands, 2000, pp.\ 307--331}
\address
{IH\'ES\\ Le Bois-Marie\\ 35, route de Chartres\\ F-91440
Bures-sur-Yvette, France\\ and Department of Mathematics\\ Michigan
State University\\ East Lansing, MI 48824-1027\\ USA}
\curraddr{School of Mathematics\\University of Minnesota\\Minneapolis,
  MN 55455\\USA}
\email{voronov@math.umn.edu}

\date{January 15, 2000; minor revision April 6, 2016}

\dedicatory{Dedicated to the memory of Mosh\'e Flato}

\begin{abstract}
  The goal of this paper is to complete Getzler-Jones' proof of
  Deligne's Conjecture, thereby establishing an explicit relationship
  between the geometry of configurations of points in the plane and
  the Hochschild complex of an associative algebra. More concretely,
  it is shown that the $B_\infty$-operad, which is generated by
  multilinear operations known to act on the Hochschild complex, is a
  quotient of a certain operad associated to the compactified
  configuration spaces. Different notions of homotopy Gerstenhaber
  algebras are discussed: One of them is a $B_\infty$-algebra,
  another, called a homotopy G-algebra, is a particular case of a
  $B_\infty$-algebra, the others, a $G_\infty$-algebra, an
  $\oE$-algebra, and a weak $G_\infty$-algebra, arise from the
  geometry of configuration spaces. Corrections to the paper of
  Kimura, Zuckerman, and the author related to the use of a nonextant
  notion of a homotopy Gerstenhaber algebra are made.
\end{abstract}

\maketitle

In an unpublished paper of E.~Getzler and J.~D.~S. Jones \cite{gj},
the notion of a homotopy $n$-algebra was introduced. Unfortunately the
construction that justified the definition contained an error, which
passed unnoticed in subsequent work, in spite of being heavily used in
it. That work included the solution by Getzler and Jones \cite{gj} of
Deligne's Conjecture, whose weak version had been proven in
\cite{gv1}; the construction by T.~Kimura, G.~Zuckerman, and the
author \cite{kvz} of a homotopy Gerstenhaber algebra structure (called
a $G_\infty$-algebra therein) on the state space of a topological
conformal field theory (TCTF); the extensions of the above work by
Akman \cite{akman:hoch,akman:master} and Gerstenhaber and the author
\cite{gv1}; a few papers delivered at the Workshop on Operads in
Osnabr\"uck in June 1998 \cite{osna}. The purpose of this paper is to
correct the error in the original construction of \cite{gj}, complete
Getzler-Jones' proof of Deligne's Conjecture accordingly, and make
appropriate corrections in \cite{kvz}.

First, let us describe the problem. A \emph{Gerstenhaber $($G-$)$
algebra} is defined by two operations, a (dot) product $ab$ and a
bracket $[a,b]$, on a graded vector space $V$ over a ground field $k$
of characteristic zero, so that the product defines a graded
commutative algebra structure on $V$ and the bracket a graded Lie
algebra structure on $V[1]$, the desuspension of the graded vector
space $V = \bigoplus_n V^n$: $V[1]^n = V^{n+1}$. The bracket must be a
graded derivation of the product in the following sense:
\[
[a , bc] = [a,b] c + (-1)^{(\abs{a} - 1) \abs{b}} b [a,c] ,
\]
where $\abs a$ denotes the degree of an element $a \in V$.  In other
words, a G-algebra is a specific graded version of a Poisson
algebra.

A G-algebra may be equivalently defined as an algebra over the operad
$e_2 = H_\bullet (D) = H_\bullet(D;k)$ of the homology of the
\emph{little disks operad}. This is a collection of manifolds $D(n)$,
$n \ge 1$, where each $D(n)$ is the configuration space of $n$
nonoverlapping little disks inside the standard unit disk in the
plane: \medskip

\centerline{\epsfxsize=1in \epsfbox{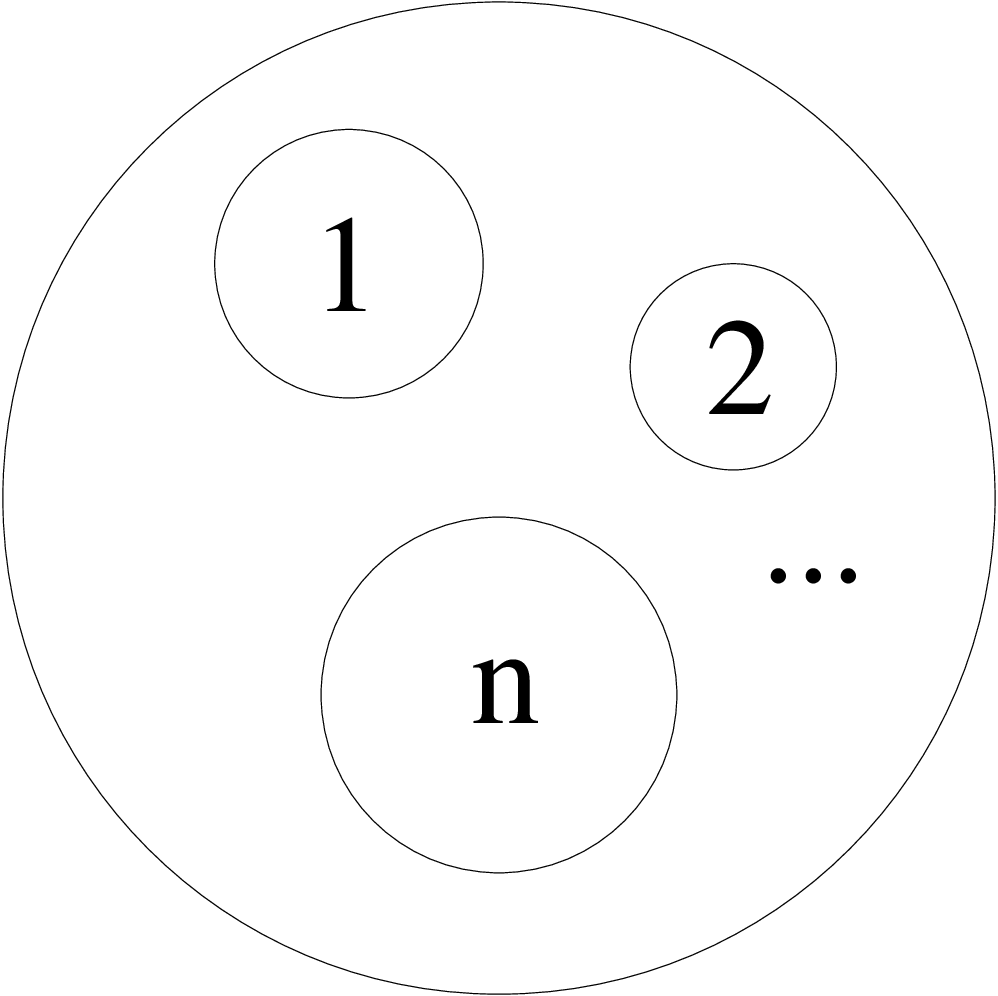}}
\smallskip

The space $D(n)$ is obviously an open set in $\nr^{3n}$ (whence a
manifold structure), each configuration being uniquely determined by
the position of the centers of the disks and their radii. It is assumed
that each little disk is labeled by a number from 1 through $n$, which
defines the action of the permutation group $\Sigma_n$ on $D(n)$. The
operad composition
\[
\gamma: D(k) \times D(n_1) \times \dots \times D(n_k) \to 
D(n_1 + \dots + n_k)
\]
is given by scaling down given configurations in $D(n_1)$, \dots,
$D(n_k)$, gluing them into the $k$ holes in a given configuration in
$D(k)$, and erasing the seams. Thus $D = \{ D(n) \; | \; n \ge 1\}$
becomes an operad of manifolds.

The way the little disks operad $D$ has relevance to G-algebras is
through the following theorem.
\begin{thm}[F. Cohen \cite{fcohen:disks}]
\label{cohen}
The structure of a G-algebra on a $\nz$-grad\-ed vector space is
equivalent to the structure of an algebra over the homology little
disks operad $H_\bullet(D)$.
\end{thm}
In view of this theorem, we will refer to the operad $e_2 :=
H_\bullet(D)$ as the \emph{G-operad}.

On the other hand, a purely algebraic example of a G-algebra was given
by the following result.
\begin{thm}[M.~Gerstenhaber \cite{gerst}]
The \h \coh\ of an associative algebra has the natural structure of a
G-algebra with respect to the cup product and the Gerstenhaber (G-)
bracket.
\end{thm}

Thus the \h \coh\ is naturally an algebra over the homology little
disks operad $e_2$. In a letter \cite{deligne:letter}, P.~Deligne
pointed out that there must be more than this formal relationship
between the little disks operad and the \h \coh, which later became
quoted as the following conjecture.
\begin{conj}[Deligne's Conjecture]
The structure of an algebra over the homology little disks operad
$e_2$ on the \h \coh\ may be naturally lifted to the $($co$)$chain
level.
\end{conj}

The goal of this paper is to prove a version of this conjecture, see
Corollary~\ref{final}. The conjecture has found several
interpretations, both algebraic and geometric. Algebraically, the
conjecture would imply the structure of a \hg algebra on the \hc\
complex. Here is an account of what is known about it up to date.

\begin{itemize}

\item
A \hg algebra structure on the \h complex was defined by Gerstenhaber
and the author \cite{gv1}. Apart from the cup product and the
G-bracket, we used higher operations (called \emph{braces}) on the \h
complex constructed by T.~Kadeishvili \cite{kadeishvili} and Getzler
\cite{g:cartan}, and wrote down a set of identities for the operations
interpreted as homotopies for the G-algebra identities and a hierarchy
of homotopies between homotopies. The bilinear and trilinear braces
providing homotopies for the commutativity of the cup product and the
distributivity of the G-bracket, respectively, had been known since
the original paper of Gerstenhaber \cite{gerst}.

\item
Getzler and Jones \cite{gj} showed that the same operations on the \h
complex defined the structure of a \emph{$B_\infty$-algebra}, a more
general version of a \hg algebra introduced by H.~Baues in the study
of the double bar construction in algebraic topology, see
Definition~\ref{baues} below. The $B_\infty$-algebra structure on the
\h complex is obtained by setting one of the $B_\infty$-operations to
be the dot product, some others to be the braces, and the others to
zero. In terms of operads, the $B_\infty$-structure on the \h complex
is obtained by realizing the homotopy G-operad of the previous
paragraph as a quotient of the $B_\infty$-operad and using the \hg
structure on the \h complex.

\item
Tamarkin \cite{tamarkin:def} extended the cup product and the G-bracket on
the Hochschild complex to the structure of a
\emph{$G_\infty$-algebra}, which is the most canonical notion of a \hg
algebra. It may be defined in terms of the corresponding operad, which
is the minimal model of the G-operad in the sense of M.~Markl
\cite{markl}. The $G_\infty$-algebra structure of Tamarkin was again
constructed by defining a morphism of operads $G_\infty \to B_\infty$
and using the $B_\infty$-algebra structure described in the previous
two paragraphs. The construction of the morphism $G_\infty \to
B_\infty$ is very involved: It uses the existence of Drinfeld's
associator and Etingof-Kazhdan's quantization theorem. Tamarkin
\cite{tamarkin:def} used this construction in his algebraic proof of
Kontsevich's Formality Theorem \cite{kon:dq}, which implied the
Deformation Quantization Conjecture \cite{bffls}.

\end{itemize}

\begin{rem}
Note that we used the same notation ``$G_\infty$-algebra'' for a
different object in our paper \cite{kvz}. Since as it has turned out,
that object does not exist, this should not create any confusion.
\end{rem}

The above results are purely algebraic. However they bring more
evidence to the relationship between the geometry of the little disks
operad and the algebra of the \h complex hinted by Deligne's
Conjecture. Moreover, to establish this relationship, one just needs
to relate the above operads (the \hg operad, $B_\infty$, and
$G_\infty$) to the little disks operad. Here is what is known in this
direction.

\begin{itemize}

\item
Getzler and Jones \cite{gj} noticed that the $G_\infty$-operad is
isomorphic to the the first term of a spectral sequence associated to
the so-called ``topological'' filtration of the moduli space operad of
configurations of points in the plane. This operad is homotopy
equivalent to the little disks operad. For the moduli space operad,
Getzler and Jones also offered the geometric construction of a
cellular model which mapped surjectively to the
$B_\infty$-operad. Unfortunately, there was an error in the
construction: The ``cellular model'' was, strictly speaking, not
cellular, because its components were not cell complexes. A
counterexample was found by Tamarkin.

\item
Another interesting relation between the $G_\infty$-operad and the
little disks operad is the idea of a recent preprint of Tamarkin
\cite{tamarkin:chains}, who showed that the singular chain operad
$C_\bullet(D)$ of little disks is \emph{formal}, \emph{i.e}.,
quasi-isomorphic to its homology $H_\bullet(D)$. Kontsevich
\cite{kon:after} found a very simple geometric proof of this
result. Both results imply that there exists a morphism $G_\infty \to
C_\bullet(D)$ unique up to homotopy.

\end{itemize}

Since all the algebraic solutions of Deligne's Conjecture so far used
the explicit $B_\infty$-structure on the Hochschild complex, finding
an explicit relationship between the $B_\infty$-operad and the little
disks operad seems to be most crucial for understanding the
relationship between algebra and geometry suggested by Deligne's
Conjecture. Moreover, the operads $G_\infty$ and $C_\bullet(D)$, which
have obvious geometric meaning, do not act explicitly on the
Hochschild complex. In our understanding, placing the
$B_\infty$-operad within the topology of the little disks operad will
provide the most complete solution of Deligne's Conjecture.

In this paper, we show that the $B_\infty$-operad is a quotient of a
certain operad $\oE$ arising geometrically from the moduli space
operad of configurations of points in the plane. The operad $\oE$ is
the first term of the spectral sequence associated to a filtration of
the moduli space operad. We show that the operad $\oE$ is free as an
operad of graded vector spaces and quasi-isomorphic to the
G-operad. We construct an explicit surjection from the operad $\oE$ to
the $B_\infty$-operad, which will imply an explicit $\oE$-algebra
structure on the Hochschild complex. We find it truly remarkable that
the object $B_\infty$ of the algebraic world (with its natural action
on the \h complex) is ruled by the geometry of the configuration
spaces (through the operad $\oE$). For example, see
Section~\ref{relation}, where the defining relations of a
$B_\infty$-algebra are read off from the incidence relations for the
strata in the moduli spaces.

Since $\oE$ is quasi-isomorphic to the G-operad, it follows from
homotopy theory of operads that there exists a morphism $G_\infty \to
\oE$. This also implies the existence of a $G_\infty$-algebra
structure on the Hochschild complex.

At the end we will indicate which changes are to be made to the paper
\cite{kvz} by Kimura, Zuckerman, and the author, which utilized the
incorrect notion of a $G_\infty$-algebra. Briefly, the changes are
that this notion should be replaced by the correct one, which brings
corrections to the identities described implicitly in
\cite{kvz}. However, all the identities written out in \cite{kvz}
explicitly do not require corrections.

\begin{rem}
While this paper was in preparation, there were made two announcements
of results of similar nature. J.~McClure and J.~Smith \cite{mcs}
constructed a cellular operad acting on the \h complex and announced
that it was homotopy equivalent to the little disks operad. Kontsevich
\cite{kon:after}, more details are coming in \cite{kon-soib},
announced the construction of an operad with an explicit
quasi-isomorphism to the G-operad and an explicit action on the
Hochschild complex, along with the proof of a multi-dimensional
generalization of Deligne's Conjecture.
\end{rem}

We will use the following \emph{terminology} regarding basic notions
of topology. All \emph{topological spaces} considered will be
Hausdorff, except when referring to complex algebraic curves, we will
use standard terminology of Zariski topology, such as an irreducible
component. For a topological space $X$, let $X^\bullet$ denote its
\emph{one-point compactification}, which is $X^\bullet = X$ if $X$ is
compact and $X^\bullet = X \cup \{\infty\}$ otherwise. A
\emph{$(p$-dimensional$)$ cell} in a topological space $X$ is a subset
$E \subset X$ along with a continuous map $f: I^p \to \bar{E}$, where
$I^p$ is the closed unit cube in $\nr^p$ with the boundary $\del I^p$,
such that $f$ is a homeomorphism in the interior of $I^p$. A
\emph{cellular partition of $X$} is a partition of $X$ into the
disjoint union of cells. A \emph{cell-complex structure on $X$}, or
equivalently, a \emph{cellular decomposition of $X$}, is a cellular
partition such that for each $p$-dimensional cell $E \subset X$, its
boundary $\del E$ is contained in $X^{p-1}$, where $X^{p-1}$ is the
\emph{$p-1$ skeleton of $X$}, the union of $p-1$-dimensional cells.
All cell complexes considered in the paper will be \emph{finite},
\emph{i.e}., consist of finitely many cells, and therefore
automatically be \emph{CW-complexes}. Cell complexes form a tensor
category with respect to cellular maps and direct products. A
\emph{cellular operad} is an operad of cell complexes. A
\emph{stratification} of a manifold is a decomposition of the manifold
into the disjoint union of connected submanifolds, called
\emph{strata}, so that the boundary of a stratum is the union of
strata of lower dimensions.

\begin{ack}
  I would like to thank F.~Akman, Z.~Fiedorowicz, Victor Ginz\-burg,
  V.~Hinich, Y.-Z. Huang, T.~Kimura, M.~Kon\-tse\-vich, M.~Markl,
  J.~Stasheff, D.~Sullivan, and D.~Tamarkin for helpful discussions. I
  am grateful to V.~Tourchine for pointing out an inaccuracy in
  defining the filtration $\overline F^\bullet$ in Section
  \ref{stratification} and N.~Bottman for noticing an error in the
  description of Tamarkin's counterexample, Remark~\ref{Tamarkin} in
  the same section, in earlier versions of this paper. I would also
  like to use this opportunity to express my sincere gratitude to
  Mosh\'e Flato for his energetic support and appreciation of the work
  of young people entering the vast grounds of Deformation
  Quantization. I also thank the Institut des Hautes \'Etudes
  Scientifiques, where the paper was finished, for hospitality.
\end{ack}

\section{The $G_\infty$-operad}

\subsection{Getzler-Jones' cellular partition}
Let us recall the construction of Getzler and Jones \cite{gj} of a
cellular partition of compactified configuration spaces. 

\subsubsection{Fox-Neuwirth cells}
\label{FN}
Consider Fox-Neuwirth's cellular decomposition \cite{fn} of the
one-point compactification of a configuration space $F(\nc,n) = \nc^n
\setminus \Delta$, where $\Delta$ is the fat diagonal $\cup_{i \ne j}
\{x_i = x_j\}$, of $n \ge 1$ distinct points in the complex plane
$\nc$: The cells, which we will call \emph{Fox-Neuwirth cells}, will
be labeled by ordered partitions of the set $\{1, \dots, n\}$ into
ordered subsets. For example, $\{3\}\{2,1\}$ denotes a partition of
$\{1,2,3\}$ into two subsets: the first subset is $\{3\}$ and the
second subset is $\{2,1\}$, ordered so that 2 precedes 1. Partitioning
labels $\{1, \dots, n\}$ into ordered subsets reflects grouping points
lying on common vertical lines on the plane. Ordering between subsets
is the left-to-right order between the vertical lines; ordering
within a subset is the bottom-to-top order within the vertical
line. For each $n \ge 1$, take the quotient space
\[
\MM (n)  = F(\nc,n)/ \nr^2 \rtimes \nr_+^*
\]
by the action of translations and dilations on $F(\nc,n)$. The
dimension of $\MM(n)$ is equal to $2n - 3$ for $n \ge 2$ and 0 for $n
= 1$. The Fox-Neuwirth cells are obviously invariant under this
action, and their quotients, which we will also call
\emph{Fox-Neuwirth cells}, make up a cellular decomposition of the
one-point compactification $\MM(n)^\bullet$ of $\MM(n)$. The spaces
$\MM(n)$ do not form an operad, but one can glue lower $\MM (k)$'s to
the boundaries of higher $\MM (n)$'s to form a topological operad $\FF
= \{\F{n}\; | \; n \ge 2\}$. In fact, the underlying spaces $\F{n}$
are smooth manifolds with corners compactifying $\MM (n)$, see next
section.

All the spaces $D(n)$, $F(\nc,n)$, $\MM(n)$, and $\F{n}$ are homotopy
equivalent. Moreover $D$ and $\F{n}$ are homotopy equivalent operads,
see \cite{gj}.

\subsubsection{Compactified moduli spaces}
The resulting space $\F{n}$ is an $S^1$-bundle over the real
compactification $\X{n+1}$ of the moduli space $\MM_{0,n+1}$ of
$n+1$-punctured curves of genus zero, see \cite{gj,gv2}. The space
$\F{n}$ can also be interpreted as a ``decorated'' moduli space, see
\cite{gv2}. Indeed, it can be identified with the moduli space of data
$(C; x_1, \dots, x_{n+1}; \linebreak[0] \tau_1, \dots, \tau_m,
\tau_\infty)$, where $C$ is a stable complex complete algebraic curve
with $n+1$ punctures $x_1, \dots, x_{n+1}$ and $m$ double points. For
each $i$, $ 1 \le i \le m$, $\tau_i$ is the choice of a tangent
direction at the $i$th double point to the irreducible component that
is farther away from the ``root'', \emph{i.e}., from the component of
$C$ containing the puncture $\infty := x_{n+1}$, while $\tau_\infty$
is a tangent direction at $\infty$. The stability of a curve is
understood in the sense of Mumford's geometric invariant theory: Each
irreducible component of $C$ must be stable, \emph{i.e}., admit no
infinitesimal automorphisms. The operad composition is given by
attaching the $\infty$ puncture on a curve to one of the other
punctures on another curve, keeping the tangent direction at each new
double point.

\subsubsection{Getzler-Jones cells}
\label{GJ}
Fox-Neuwirth's cell decomposition of each space $\MM (n)^\bullet$
induces a cell partition of the compactification $\F{n}$ in the
following way. First of all, by an $n$-\emph{tree} we mean a directed
rooted tree with $n$ labeled initial edges, the \emph{leaves}, and one
terminal edge, the \emph{root}, each of these $n+1$ edges incident to
only one vertex of the tree, such that the number $n (v)$ of the
incoming edges for any vertex $v$ is at least two. Define the
\emph{tree degree} of an $n$-tree $T$ as $n-v(T)-1$, where $v(T)$ is
the number of vertices in $T$. The cells, which we will call
\emph{Getzler-Jones cells}, in $\F{n}$ are enumerated by pairs $(T,
p)$, where $T$ is a tree, labeling a stratum of the boundary of
$\F{n}$, and $p$ is a function $p(v)$ on the vertices $v$ of the tree
$T$, such that each $p(v)$ is an ordered partition, as in \ref{FN}
above, of the set $\IN (v)$ of incoming edges for a vertex $v$ of the
$T$. These partitions $p(v)$ label cells in the corresponding open
moduli spaces $\MM (n (v))$, whose products make up the
stratum. However generally speaking, it is not true that this cellular
partition is a cell complex: The boundary of a $q$-dimensional cell
does not always lie in the $q-1$-skeleton. We will take the union of
certain Getzler-Jones cells to form a stratification of $\F{n}$
compatible with the operad structure.

\subsection{Stratification of $\FF$}

\subsubsection{Stratification of $\MM$}
\label{K}

Consider the following subsets in $\MM(n)$.
\begin{itemize}
\item
For each ordering of the set $\{1, \dots, n\}$, consider the
corresponding Fox-Neu\-wirth cell. It consists of configurations of $n$
points on a vertical line in the prescribed order going from bottom to
top.
\item
For each ordered partition of the set $\{1, \dots, n\}$ into two
parts, consider the corresponding Fox-Neu\-wirth cell. It consists of
configurations of $n$ points on two vertical lines in the prescribed
order going from bottom to top and from left to right.
\item
The complement to the union of subsets of the above two types.
\end{itemize}
These subsets obviously form a stratification of the manifold
$\MM(n)$. We will also need the following filtration of $\MM(n)$ into
three closed subsets: The closure $J_1$ of the union of the strata of
the first type, the closure $J_2$ of the union of the strata of the
second type, and $J_3$ which is all of $\MM(n)$.

Adding the point at $\infty$ to each filtration component and setting
$J_0 = \{\infty\}$, we get a filtration of the pointed space
$\MM(n)^\bullet$ for $n \ge 3$:
\begin{equation}
\label{filtr}
J_0 \subset J_1 \subset J_2 \subset J_3 = \MM (n)^\bullet .
\end{equation}
For $n=2$, the space $\MM(2)$ is diffeomorphic to $S^1$ and thereby
compact. In this case we will not add a point $\infty$ to any
filtration components and will set $J_0 = \emptyset$.

\begin{sloppypar}
Consider the corresponding homological spectral sequence converging to
$H_\bullet (\MM(n)^\bullet, \infty)$, which is naturally isomorphic to
$H^\bullet (\MM(n))$ by Poincar\'e-Lefschetz duality, all coefficients
being taken in $k$. The term $E^1$ may be identified with the sum of
\[
E^1_{p,q} = H_{p+q}(J_p, J_{p-1}), \qquad 1 \le p \le 3, \, -p \le q
\le 2n -p - 3.
\]
\end{sloppypar}

\begin{prop}
\label{collapse}
The homological spectral sequence associated to the filtration
\eqref{filtr} collapses at $E^2$.
\end{prop}

\begin{proof}
Since $J_2 \setminus J_1$ and $J_1 \setminus J_0$ are disjoint unions
of $n! (n-1)$ and $n!$ cells of dimension $n-1$ and $n-2$,
respectively,
\[
\dim H_{2+q} (J_2, J_1) = \begin{cases}
        n! (n-1) & \text{for $2+q = n-1$},\\
        0 & \text{otherwise},
\end{cases}
\]
and
\[
\dim H_{1+q} (J_1, J_0) = \begin{cases}
        n!  & \text{for $1+q = n-2$},\\
        0 & \text{otherwise}.
\end{cases}
\]
Note also that with respect to the Fox-Neuwirth cellular decomposition
of $J_3$, the subspace $J_2$ is a cell subcomplex. The dimensions of
cells in the complement of $J_2$ run from $n$ through $2n-3$. Thus,
\[
H_{p+q} (J_3, J_2) =    0 \qquad \text{unless $n \le p+q \le 2n-3$}.
\]
Observe that for any pair $(p,q)$ and integer $r \ge 2$, either
$E^r_{p,q}$ or $E^r_{p-r,q+r-1}$ is 0. Thus the spectral-sequence
differentials $d^r: E^r_{p,q} \to E^r_{p-r,q+r-1}$ will all vanish for
$r \ge 2$, which implies the collapse of the spectral sequence at
$E^2$.
\end{proof}

\subsubsection{Stratification of $\FF$}
\label{stratification}

First of all, stratify $\F{n}$ by the topological type of the stable
algebraic curve. We will refer to this stratification as
\emph{coarse}. Each stratum $S_T$ will correspond to a tree $T$ and be
isomorphic to the product of the spaces $\MM (n (v))$ over the
vertices $v$ of the tree as in \ref{GJ}. Then subdivide the coarse
stratification as follows. For each space $\MM(n (v))$, take the
stratification into three types of strata as in the previous
section. The products of these strata over the set of vertices $v$ of
$T$ will form a finer partition of the space $\F{n}$.  Below is a
figure denoting the part obtained as the product of two strata $J_2
\setminus J_1$ in $\MM(3)$ and $J_1 \setminus J_0$ in $\MM(2)$.

\noindent
\let\picnaturalsize=N
\def\picsize{3.5in}
\def\picfilename{sphere.eps}
\ifx\nopictures Y\else{\ifx\epsfloaded Y\else\input epsf \fi
\let\epsfloaded=Y
\centerline{\ifx\picnaturalsize N\epsfxsize \picsize\fi
\epsfbox{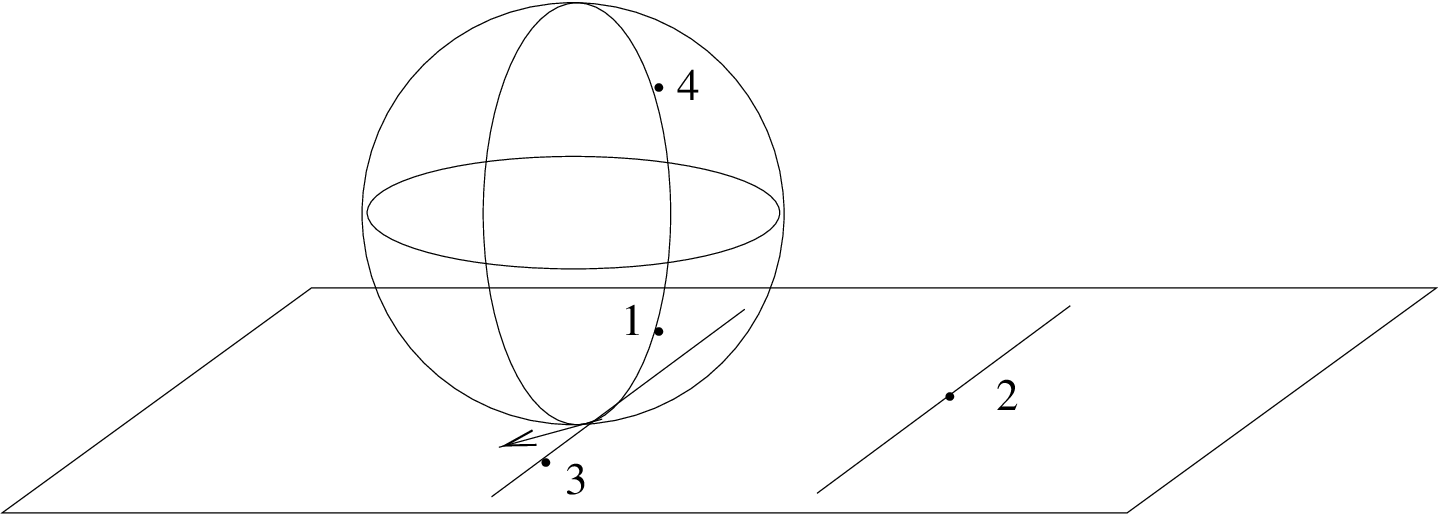}}}\fi

\noindent
Here the plane denotes the irreducible component of the stable curve
with the puncture $\infty$ placed at $\infty$. The other irreducible
component is the sphere. The arrow at the double point determines the
positive direction of the real axis on the sphere. We prefer to work
with the following replacement of the above figure, which one may
think of as projection of the three-dimensional figure onto the plane
after rotating the sphere with the arrow, so that the arrow points in
the direction of the positive real axis on the plane. The circle may
be thought of as a magnifying glass through which the observer living
on the plane sees what happens in the infinitesimal world at the
double point.

\noindent
\let\picnaturalsize=N
\def\picsize{1.5in}
\def\picfilename{proj.eps}
\ifx\nopictures Y\else{\ifx\epsfloaded Y\else\input epsf \fi
\let\epsfloaded=Y
\centerline{\ifx\picnaturalsize N\epsfxsize \picsize\fi
\epsfbox{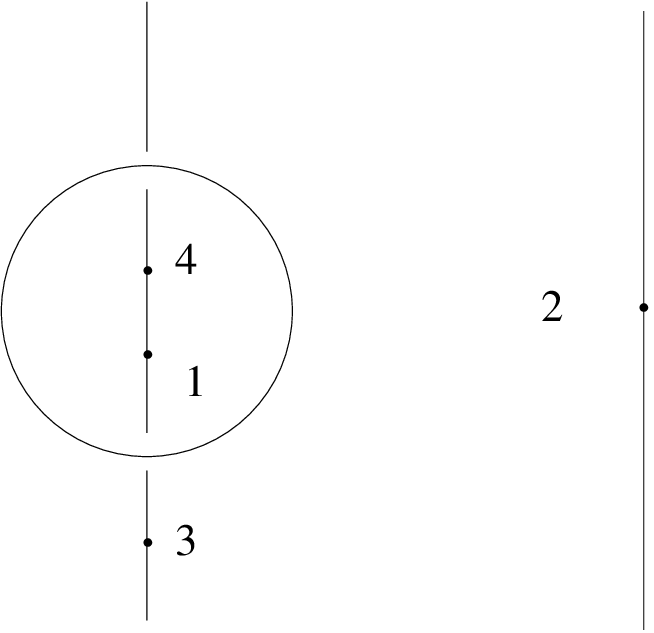}}}\fi

\noindent
We will usually think of either of these figures as not a single
configuration of points, but rather the set of all configurations of
points which are in the relative position determined by the figure.

The above partition of the space $\F{n}$ will in fact be a
stratification. The reason is that the boundary of each part may be
split into the union of the boundary within the coarse stratum $S_T$
and the boundary within the boundary of $S_T$. The boundary of the
first type is the union of other strata within the coarse stratum
$S_T$, because we used a stratification of each $\MM(n(v))$ to
construct our partition of $\F{n}$. The boundary of the second type is
obtained by letting a few groups of points on the irreducible
components of the stable curve bubble off, forming new components
attached at double points. If all the points on a component lie on a
single vertical line, the points on the limiting components will also
group on a single vertical line within each bubble, and again, the
whole fine part will be in the boundary. If all the points on a
component lie on two vertical lines, then the points on the components
that bubble off will lie on two or one vertical lines, and again, the
whole fine part will be in the boundary, because by translations and
dilations, one can always match up the vertical lines on different
components (if less than three on each), for example, pass through
points 0 or 1 on the real axis. Since $J_3\setminus J_2$ comprises the
top strata in $\MM(n(v))$, the closure of this component of the
stratum will be all $\F{n(v)}$, therefore, the boundary is $J_2 \cup
\del \F{n(v)}$, which is the union of lower-dimensional strata.

\begin{rem}
\label{Tamarkin}
  Note that if we extended our partition to a cell partition by taking
  all the Getzler-Jones cells, the boundary of a cell would not be the
  union of other cells, in general. For example, take the Fox-Neuwirth
  cell formed by six points on three vertical lines, two points on
  each. Its boundary has a nonempty intersection with the following
  Getzler-Jones cell:

\noindent
\let\picnaturalsize=N
\def\picsize{.7in}
\def\picfilename{counter.eps}
\ifx\nopictures Y\else{\ifx\epsfloaded Y\else\input epsf \fi
\let\epsfloaded=Y
\centerline{\ifx\picnaturalsize N\epsfxsize \picsize\fi
\epsfbox{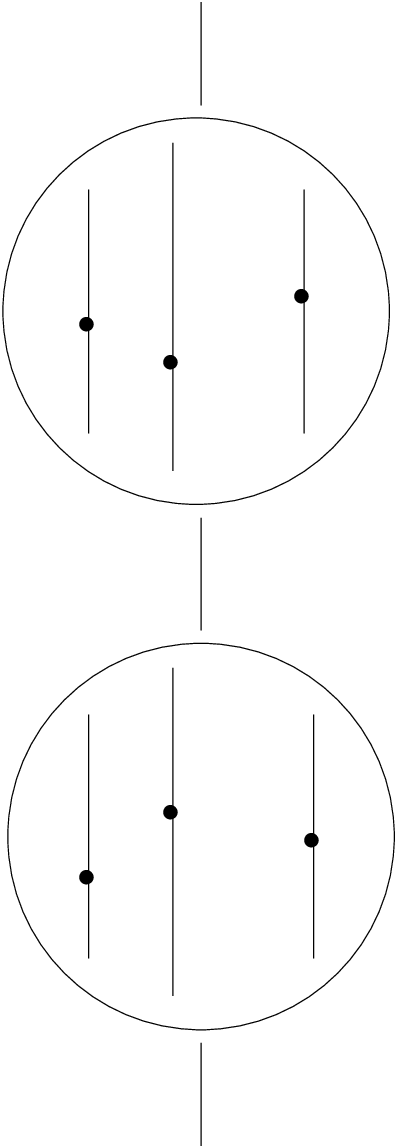}}}\fi

\noindent

However, only part of the this Getzler-Jones cell will be in the
boundary of the Fox-Neuwirth cell. This part will consist of those
positions of three vertical lines on the two bubbles that can be
brought together by the action of $\nr^2 \rtimes \nr_+^*$ on a single
copy of the Riemann sphere. This is exactly the problem with
Getzler-Jones' partition not being a cell complex. Our stratification
is designed to get around this problem.
\end{rem}

Now form two filtrations $\dots \subset F^p \subset F^{p+1} \subset
\dots $ and $\dots \subset \overline F^p \subset \overline F^{p+1}
\subset \dots $. The \emph{topological} filtration $F^\bullet$ is
obtained by taking $F^p$ to be the closure of the union of coarse
strata of dimension $p$. The filtration $\overline F^\bullet$ is
defined as follows:
\[
\overline F^p = \coprod_{\text{$n$-trees T}} \coprod_{i(v)} \prod_{v
\in T} J_{i(v)} (n(v)),
\]
where the union is over all functions $i(v)$ from the set of vertices
$v$ of the tree $T$ to the set $\{1,2,3\}$, satisfying the condition
$\sum_{v \in T} (n(v) + i(v) - 3) \le p$. Note that when $i(v) = 1$ or
$2$, $n(v) + i(v) - 3 = \dim J_{i(v)} (n(v))$.

Consider the corresponding spectral sequences, which converge to
$H_\bullet (\F{n})$ and whose first terms are
\begin{align*}
E^1_{p,q} = H_{p+q} (F^p, F^{p-1}),\\
\overline E^1_{p,q} = H_{p+q} (\overline F^p, \overline F^{p-1}) .
\end{align*}

\subsubsection{Operad properties of the first spectral sequence $E^r$}

Note that the operad composition respects the two filtrations of
spaces $\F{n}$, $n \ge 1$. Therefore, we are getting two operads of
spectral sequences, cf. \cite{ksv2}. In particular, the first terms
$E^1$ and $\overline E^1$ of the spectral sequences form operads of
complexes. There is a purely algebraic interpretation of the first
operad $E^1$, noticed by Getzler and Jones in \cite{gj}.

\begin{prop}[Getzler and Jones]
\label{prop1.2}
The operad $E^1$ is naturally isomorphic to the cobar construction of
the Gerstenhaber operad $e_2 = H_\bullet (D)$.
\end{prop}

\begin{proof}
Notice that for the space $\F{n}$, $n \ge 2$, the group $E^1_{p,q} =
H_{p+q} (F^p, \linebreak[0] F^{p-1})$ is naturally isomorphic to
$H^{-q} (F^p \setminus F^{p-1})$ by Poincar\'e-Lefschetz duality. The
space $F^p \setminus F^{p-1}$ is the disjoint union of the strata
$S_T$, $T$ running over the set of $n$-trees $T$ which have exactly
$2n-2-p$ vertices. Each stratum $S_T$ is naturally isomorphic to the
product of spaces $\MM(n(v))$ over the set of vertices $v$ of the tree
$T$. Each space $\MM(n(v))$ is homotopy equivalent to the
configuration space $D(n(v))$ of $n(v)$ little disks. Thus
\begin{equation}
\label{E^1}
E^1_{p,q} = \bigoplus_{\text{$n$-trees $T$}} \left( \bigotimes_{v \in
T} H^\bullet (D(n(v))) \right)^{-q},
\end{equation}
where the superscript $-q$ means the component of degree $-q$. The
differential
\[
d^1: E^1_{p,q} \to E^1_{p-1,q}
\]
takes the component $\bigotimes_{v \in T} H^\bullet (D(n(v)))$
corresponding to a $n$-tree $T$ to the sum of components
$\bigotimes_{v \in \widehat T} H^\bullet (D(n(v)))$ over all $n$-trees
$\widehat T$ such that the tree $T$ may be obtained by contracting an
interior edge of $\widehat T$, merging two adjacent vertices $v_1$ and
$v_2$ on the tree $\widehat T$ into a single vertex $v_3$ of $T$. The
matrix element $d^1_{T \widehat T}$ of the differential is induced (up
to a sign, which is treated below) by the map
\[
\circ_i^*: H^\bullet (D(n(v_3))) \to H^\bullet (D(n(v_1))) \otimes H^\bullet
(D(n(v_2)))
\]
which is the dual of the corresponding operad structure map
\[
\circ_{i}: H_\bullet (D(n(v_1))) \otimes H_\bullet (D(n(v_2))) \to
 H_\bullet (D(n(v_3))).
\]
This map is induced on homology by the map
\[
\circ_{i}:D(n(v_1)) \times D(n(v_2)) \to  D(n(v_3))
\]
gluing the unit disk with a configuration of $n(v_2)$ little disks
into the $i$th little disk in a configuration of $n(v_2)$ little
disks, $i$ corresponding to the contracted edge in $\IN (v_1)$, if we
assume that the contracted edge is directed from $v_2$ to $v_1$.

The sign for $d^1_{T \widehat T}$ comes from the choice of orientation
of the strata $S_T$. This orientation may be chosen by ordering all
the edges of each $n$-tree $T$, except the root edge. The orientation
on the stratum $S_T$ is then given by ordering the $x$ and the $y$
coordinates of the points in $\nc$ corresponding to the edges of $T$,
according to the order of the edges, skipping for each vertex $v$ the
coordinates of the point corresponding to the first edge and the $x$
coordinate of the point corresponding to the second edge --- remember
that $\MM (n) = F(\nc,n) / \nr^2 \rtimes \nr_+^*$. Then the
compatibility of orientations on $S_T$ and $S_{\widehat T}$ implies
that $d_{T \widehat T}^1$ is $\circ_i^*$ multiplied by the sign of the
permutation from the ordered set of edges of $\widehat T$ to the
ordered set $\{e, \text{edges of $T$}\}$, where $e$ is the contracted
edge in $\widehat T$.

This description of $E^1$ in terms of trees and the operad $H_\bullet
(D)$ of vector spaces means that $E^1$ is the cobar construction of
the operad $H_\bullet (D)$, just by the definition of Ginzburg-Kapranov
\cite{gk}.

\end{proof}

\begin{crl}[Getzler and Jones, Markl]
\label{crl}
The operad $E^1$ is a free resolution of the G-operad $e_2$,
\emph{i.e}., there is a morphism of operads
\[
E^1 \to e_2
\]
inducing an isomorphism on homology and $E^1$ is free as an operad of
graded vector spaces. Moreover $E^1$ is a minimal model of $e_2$.
\end{crl}

\begin{rem}
Here we sketch a proof due to Markl \cite{markl:dist}, which is
different from that of Getzler and Jones \cite{gj}. The fact that
$E^1$ is a minimal model of $e_2$ was first noticed by Markl
\cite{markl:dist}.
\end{rem}

\begin{proof}
Using the description of the operad $e_2$ in terms of generators and
relations, it is a straightforward exercise to check that the
quadratic dual of $e_2$ is again $e_2$ up to a shift of grading and
the change of it to the opposite. Then from the proof of
Proposition~\ref{prop1.2}, one can identify $E^1$ with the cobar
construction of $e_2^!$. The natural homomorphism, see Lemma 4.1.2 of
\cite{gk}, which easily generalizes to the graded case, from the cobar
construction of an operad $\PP$ to the quadratic dual $\PP^!$ gives
for $\PP = e_2^!$ a morphism of operads
\begin{equation}
\label{koszul}
E^1 \to e_2.
\end{equation}
It is also known that the operad $e_2$ is Koszul, see \cite{gj} or a
purely algebraic proof by Markl \cite{markl:dist}. This means (by
definition of \cite{gk}) that the morphism \eqref{koszul} is a
quasi-isomorphism.

According to Markl \cite{markl}, the cobar construction of the
quadratic dual of a Koszul operad (see \cite{gk}) is a minimal model
of that operad. Applied to $e_2$, this implies that $E^1$ is a minimal
model of $e_2$.
\end{proof}

\begin{df}
\label{1.1}
The $G_\infty$-\emph{operad} is the operad $E^1$. A
$G_\infty$-\emph{algebra} is an algebra over the $G_\infty$-operad.
\end{df}

\subsubsection{Operad properties of the second spectral sequence
$\overline E^r$}

\begin{thm}
\label{oE}
The operad $\oE$ is free as an operad of graded vector spaces, and its
homology is isomorphic to $e_2$.
\end{thm}

\begin{proof}
1. First of all, let us prove that $\oE$ is free. Recall that
\[
\oE_{p,q} = H_{p+q}(\overline F^p, \overline F^{p-1}) = H^{\dim
\overline F^p - p - q} (\overline F^p \setminus \overline F^{p-1}).
\]
Since $\overline F^{p}$'s are made out of a stratification of $\F{n}$
(see Section \ref{stratification}), we have
\[
\coprod_p \overline F^p \setminus \overline F^{p-1} = 
\coprod_{\text{$n$-trees T}} \prod_{v \in T} \coprod_{i=1}^3
J_i (n(v)) \setminus J_{i-1} (n(v)),
\]
where for each $l \ge 2$,
\begin{equation}
\label{filtrK}
J_0 (l) \subset J_1 (l) \subset J_2 (l) \subset J_3 (l)  = \MM (l)^\bullet
\end{equation}
is the filtration from Section~\ref{K}. Therefore, passing to
cohomology, we have
\[
\oE_{\bullet,\bullet} = \bigoplus_{\text{$n$-trees T}} \bigotimes_{v \in T}
\bigoplus_{i=1}^3
H^\bullet(J_i (n(v)) \setminus J_{i-1} (n(v)) ),
\]
which by definition means that $\oE_{\bullet,\bullet}$ is a free
operad generated by the collection
\begin{equation}
\label{generators}
\bigoplus_{i=1}^3
H^\bullet(J_i (n) \setminus J_{i-1} (n) ) = \bigoplus_{i=1}^3
H_\bullet(J_i (n), J_{i-1} (n) )
\end{equation}
of graded vector spaces with an action of the symmetric group $S_n$, $
n \ge 2$.

2. The next step is to show that the spectral sequence $\overline E^r$
collapses at the second term $\overline E^2$. In order to show that
the cohomology of $\oE$ is $\overline E^\infty = e_2$, regard $\oE$ as
a filtered complex, with the $k$th filtration component defined by the
tree degree $n - v(T) - 1 \le k$. We will compute the homology of
$\oE$ using the spectral sequence associated with this filtration. The
first term of this spectral sequence is $\bigoplus_{\text{$n$-trees
T}} \bigotimes_{v \in T} H^\bullet(\MM (n(v)) )$ with the Gysin
homomorphism as the differential, because of
Proposition~\ref{collapse}, that is the first term $E^1$ of the
spectral sequence associated with the filtration $F^\bullet$, see
\eqref{E^1}. Corollary~\ref{crl} shows that the homology of $E^1$ is
isomorphic to $e_2$. By construction this is the second term of the
spectral sequence associated to the filtered complex $\oE$ and the
spectral sequence converges to the homology $\overline E^2$ of
$d^1$. On the other hand, $e_2$ is the $\infty$ term $\overline
E^\infty$. Thus $e_2$ is the second term of a spectral sequence
converging to the second term of another spectral sequence converging
to the same $e_2$. This implies that both spectral sequences collapse
at the second terms. Therefore, the homology of $\oE$ is $e_2$.
\end{proof}

\section{The $B_\infty$-operad}

\subsection{$B_\infty$-algebras and the $B_\infty$-operad}

Let $V = \bigoplus_{n \in \nz} V^n$ be a graded vector space over a
field $k$ of characteristic zero, $V[1]$ its \emph{desuspension}:
$V[1] = \bigoplus_{n \in
\nz} V[1]^n$, where $V[1]^n = V^{n+1}$, and $TV[1] = \sum_{p=0}^\infty
(V[1])^{\otimes p}$ the \emph{tensor coalgebra} on $V[1]$. We will
adopt the standard notation $[a_1|\dots | a_p]$ for an element $a_1
\otimes \dots \otimes a_p \in (V[1])^{\otimes p} \subset TV[1]$. By
definition $\abs{[a_1|\dots | a_p]} = \abs{a_1} + \dots + \abs{a_p} -
p$, where $\abs a$ denotes the degree of an element $a$ in a graded
vector space. The graded coalgebra structure on $TV[1]$ is given by
the coproduct $\Delta: TV[1] \to TV[1] \otimes TV[1]$,
\[
\Delta [a_1 | \dots | a_p] = \sum_{i=0}^p [a_1 | \dots | a_i]
\otimes [a_{i+1} | \dots | a_p],
\]
for which the natural augmentation $TV[1] \to k$ is a counit.

We will be interested in studying a certain DG bialgebra structure on
$TV[1]$. Here a \emph{DG bialgebra} is an algebra $A$ with a unit, the
structure of a coalgebra, and a differential $D: A \to A[1]$, $D^2 =
0$, such that $D$ is a graded derivation and coderivation and the
comultiplication $A \to A \otimes A$ is a morphism of algebras.

\begin{df}[H. J. Baues \cite{baues}]
\label{baues}
A \emph{$B_\infty$-algebra} structure on a graded vector space $V$ is
the structure of a DG bialgebra on the tensor coalgebra $TV[1]$, such
that the element $[\,] \in (V[1])^{\otimes 0} \subset TV[1]$ is a unit
element.
\end{df}

Since the tensor coalgebra is cofree and both the differential $D:
TV[1] \to TV[1]$ and the product $M: TV[1] \otimes TV[1] \to TV[1]$
are respect the coproduct, they are determined by the compositions
\[
\pr D = \sum_{k=0}^\infty M_k: TV[1] \to V[2]
\]
and
\[
\pr M = \sum_{k,l=0}^\infty M_{k,l}: TV[1] \otimes TV[1] \to V[1]
\]
with the natural projection $\pr: TV[1] \to V[1]$. The condition $D^2
= 0$ can be rewritten as a collection of identities for the operations
$M_k$, the associativity condition for $M$ and the unit axiom for
$[\,]$ as a collection of identities for the operations $M_{k,l}$ and
the derivation property for $D$ with respect to $M$ as a collection of
identities between $M_{k,l}$ and ${M_k}$. The restriction $M_0$ of $D$
to $V[1]^{\otimes 0}$ must vanish, because $D$ must annihilate the
unit $[\,]$. The equation $D^2=0$ then implies $M_1^2 = 0$, which
means $d= M_1$ must be a differential on the graded vector space $V$,
defining the structure of a complex with a differential of degree
1. Thus, a $B_\infty$-algebra structure on a graded vector space is
equivalent to a differential $d$ and a collection of multilinear
operations $M_k$ of degree $\abs{M_k} = 2-k$ and $M_{k,l}$ of degree
$\abs{M_{k,l}} = 1-k-l$ satisfying certain identities, \emph{i.e}.,
the structure of an algebra over the DG operad generated by $M_k$, $k
\ge 2$, and $M_{k,l}$, $k,l \ge 0$, with those identities being the
defining relations. We will call this operad the
\emph{$B_\infty$-operad}.

\subsection{The algebraic description of the $B_\infty$-operad}

Here we will describe the $B_\infty$-operad explicitly. We will use
this description in the next section to show that the
$B_\infty$-operad is a quotient of the operad $\oE$, associated to the
little disks operad. Just to make the formulas more transparent, we
will describe the identities satisfied by the operations $M_k$, $k \ge
0$, and $M_{k,l}$, $k,l \ge 0$, in a $B_\infty$-algebra $V$. As we
already noticed, $M_0 = 0$ and $M_1$ is a differential $d$ on $V$. We
will adopt the following convention:
\[
M_k(a_1, \dots, a_k) := (-1)^{(k-1)\abs{a_1} + (k-2)\abs{a_2} + \dots +
\abs{a_{k-1}}} M_k [a_1| \dots | a_k],
\]
which morally means that the vertical bar $|$ has degree one and on
the left-hand side all the bars are moved between $M_k$ and
$a_1$. Here $|a_i|$ denotes the degree of $a_i$ in $V$.
However, we set
\[
M_{k,l} (a_1, \dots, a_k;
b_1, \dots, b_l) := 
M_{k,l} ([a_1| \dots | a_k] \otimes [b_1| \dots | b_l]).
\]

\subsubsection{$D^2 = 0$}

The condition $D^2 = 0$ is equivalent for the operations $M_k$, $k\ge
1$, to define an $A_\infty$-structure:
\begin{multline}
\label{A_infty}
\sum_{i+j = n+1} \sum_{k=0}^{n-j} (-1)^{\epsilon} 
M_i (a_1, \dots ,a_k, \\
M_j (a_{k+1}, \dots, a_{k+j}), \dots, a_n) = 0, \quad n \ge 1,
\end{multline}
where $\epsilon = (i+1)j + (j+1)k + i|a_1|+ (i-1)\abs{a_2} + \dots +
(i-k+1)\abs{a_{k}} + (n-k-1) |a_{k+1}| + (n- k -2) \abs{a_{k+2}} +
\dots + \abs{a_{n-1}}$ and $a_1, \dots , a_n \in V$. In fact, the sign
$\epsilon$ is obtained as $\abs{a_1} + \dots + \abs {a_k} - k$ plus
the sign coming from moving the vertical bars in any occurrence of
$M_p[a_1 | \dots | a_p]$ to the place between $M_p$ and $a_1$,
thinking of a bar as having degree 1.

\subsubsection{$[\,]$ is a unit for $M$}

This is equivalent to $M_{1,0} = M_{0,1} = \id$, $M_{k,0} = M_{0,k} =
0$ for $k \ne 1$.

\subsubsection{The associativity of $M$}

\begin{sloppypar}
The associativity of $M = \sum_{k,l \ge 0} M_{k,l}$ is equivalent to
the following identities
\begin{multline}
\label{assoc}
\begin{aligned}
\sum_{r =1}^{l+m} \sum_{\substack{l_1 + \dots + l_r = l\\ m_1 + \dots
+ m_r = m}} (-1)^\epsilon M_{k,r} ( a_1, \dots, a_k; M_{l_1,m_1}(b_1, \dots,
b_{l_1}; c_1, \dots, c_{m_1}), \\ \dots , M_{l_r,m_r}(b_{l_1 + \dots
+l_{r-1} +1}, \dots, b_{l};c_{m_1 + \dots + m_{r-1} +1}, \dots, c_m))
\end{aligned}
\\
\begin{aligned}
= \sum_{s = 1}^{k+l} \sum_{\substack{k_1 + \dots + k_s = k\\l_1 + \dots +
l_s = l}} (-1)^\delta M_{s,m} (M_{k_1,l_1}(a_1, \dots, a_{k_1};b_1, \dots,
b_{l_1}), \dots , \\
M_{k_s,l_s}(a_{k_1 + \dots + k_{s-1} +1}, \dots,
a_{k};b_{l_1 + \dots + l_{s-1} +1}, \dots, b_l); c_1, \dots, c_m)
\end{aligned}
\end{multline}
for $a_1, \dots, a_k$, $b_1, \dots, b_l$, and $c_1, \dots, c_m$ in
$V$. The sign $(-1)^\epsilon$ is the sign picked up by reordering
$[a_1 | \dots |a_k|b_1| \dots |b_l| c_1 | \dots | c_m]$ into $[a_1|
\dots |a_k |\linebreak [0] b_1 | \dots | b_{l_1} | c_1 | \dots |
c_{m_1} |\dots \linebreak[0] | b_{l_1 + \dots +l_{r-1} +1} | \dots |
b_{l} | c_{m_1 + \dots + m_{r-1} +1} | \dots | c_m]$ in the graded
vector space $TV[1]$. Similarly, $(-1)^\delta$ is the sign of
reordering $[a_1 | \dots |a_k|b_1| \dots |b_l| c_1 | \dots | c_m]$
into $[a_1 | \dots| a_{k_1}|b_1| \dots | b_{l_1} | \dots \linebreak[1]
| a_{k_1 + \dots + k_{s-1} +1}| \linebreak[0] \dots | a_{k} | b_{l_1 +
\dots + l_{s-1} +1} | \dots | b_l | c_1 | \dots | c_m]$.
\end{sloppypar}

\subsubsection{The Leibniz rule for $D$ with respect to $M$}

The fact that $D$ is a derivation of the product $M$ on $TV[1]$ is
equivalent to the following identities
\begin{multline}
\label{Leibniz}
\begin{aligned}
\sum_{n = 1}^{k+l} \sum_{\substack{k_1 + \dots + k_n = k\\l_1 + \dots
+ l_n =l}} (-1)^\epsilon M_n (M_{k_1,l_1}(a_1, \dots, a_{k_1}; b_1, \dots, b_{l_1}),
\dots , \\ M_{k_n,l_n}(a_{k_1 + \dots + k_{n-1} +1}, \dots, a_{k};
b_{l_1 + \dots + l_{n-1} +1}, \dots, b_{l}))
\end{aligned}\\
\begin{aligned}
= \sum_{r=1}^k
\sum_{i=0}^{k-r} (-1)^{\delta}
M_{k-r+1,l}(a_1, \dots, a_i, M_r (a_{i+1}, \dots, a_{i+r}), \dots, a_k;
\\
b_1 , \dots, b_l)
\end{aligned}
\\
\begin{aligned}
+ (-1)^{|a_1| + \dots |a_k| - k} \sum_{s=1}^l
\sum_{i=0}^{l-s} (-1)^{\eta}
M_{k,l-s+1}(a_1, \dots, a_k; b_1, \dots, b_i, \\
M_s (b_{i+1}, \dots,
b_{i+s}), \dots, b_l)
\end{aligned}
\end{multline}
for $a_1, \dots, a_k$ and $b_1, \dots, b_l$ in $V$. The sign
$(-1)^\epsilon$ is the sign of reordering $[a_1 | \dots
|a_k|\linebreak[0] b_1| \dots |b_l]$ into $[a_1 | \dots| a_{k_1}|b_1|
\dots | b_{l_1} | \dots | a_{k_1 + \dots + k_{n-1} +1}| \dots
\linebreak[0] | a_{k} | b_{l_1 + \dots + l_{n-1} +1} | \dots | b_l]$
in $TV[1]$ multiplied by the sign of moving $n-1$ bars between
$M_{k_1,l_1}(\dots), \dots ,\linebreak[0] M_{k_n,l_n} (\dots)$ to the
place between $M_n$ and $M_{k_1,l_1}$.  The sign $(-1)^\delta$ is
equal to $\abs{a_1} + \dots + \abs {a_i} - i$ plus the sign coming
from moving the vertical bars in $M_r[a_{i+1} | \dots | a_{i+r}]$ to
the place between $M_r$ and $a_{i+1}$. Similarly, the sign $(-1)^\eta$
is equal to $|b_1| + \dots |b_i| - i$ plus the sign coming from moving
the vertical bars in $M_s[b_{i+1} | \dots | b_{i+s}]$ to the place
between $M_s$ and $b_{i+1}$.

\begin{rem}
A few ``lower'' identities including the derivation property of the
differential $d=M_1$ with respect to the ``dot'' product $M_2$ and
identities providing homotopies for such classical identities for
binary operations as the commutativity and the associativity of the
``dot'' product, the homotopy left and right Leibniz rules for the
``circle'' product $M_{1,1}$ with respect to the ``dot'' product, and
the homotopy Jacobi identity for the ``bracket'' $[a,b] = M_{1,1}
(a,b) - (-1)^{(\abs a - 1) (\abs b - 1)} M_{1,1} (b,a)$ were written
out explicitly in \cite[Section~4.2]{kvz}. Strictly speaking, those
identities were claimed to be identities for another type of algebra,
which later turned out to be nonextant. However, one can see from
Theorem~\ref{thm:relation} of the next section, that those identities
are satisfied in a $B_\infty$-algebra.
\end{rem}

\subsection{Relation between the $B_\infty$-operad and $\oE$}
\label{relation}

\begin{sloppypar}
  The algebraic description above of the $B_\infty$-operad might be a
  good exercise in tensor algebra, but is far from inspiring. However,
  everything falls into its place, when geometry comes into play.
  Since $\oE$ is a homological operad (the degree of the differential
  is $-1$), if $B_\infty$ is a cohomological operad (the degree of the
  differential is $+1$), let us change the grading on $\oE$ to the
  opposite one, an element of degree $k$ will be assigned degree $-k$,
  from now on, so that the differential on $\oE$ is of degree $+1$.
\end{sloppypar}

\begin{thm}
\label{thm:relation}
There exists a surjective morphism of DG operads $\oE \linebreak[0]
\to B_\infty$.
\end{thm}

\begin{proof}

The operad $\oE$ is freely generated by the spaces $\bigoplus_{i=1}^3
H_\bullet(J_i (n), \linebreak[0] J_{i-1} (n) )$, $n \ge 2$, see
\eqref{generators}. Thus to define a morphism $\oE \to B_\infty$, it
suffices to define maps
\begin{equation}
\label{define}
\bigoplus_{i=1}^3
H_\bullet(J_i (n), J_{i-1} (n) ) \to B_\infty(n), \qquad n \ge2,
\end{equation}
respecting the gradings, the symmetric group actions, and the
differentials.

The complements $J_1\setminus J_0$ and $J_2 \setminus J_1$ are
disjoint unions of Fox-Neuwirth cells, and the spaces $H_\bullet (J_1,
J_0)$ and $H_\bullet (J_2, J_1)$ have the Fox-Neuwirth cells of the
types $\{i_1, \dots, i_n\}$ and $\{i_1, \dots, i_{p}\}\{i_{p+1},
\dots, i_n\}$, respectively, see \ref{FN}, as natural bases. Define the maps
\eqref{define} as follows:
\begin{align}
\label{corr1}
\{1, 2, \dots, n\} & \mapsto M_n \qquad \text{for $n \ge 2$},
\\
\label{corr2}
\{1,2, \dots, k\}\{k+1, \dots, k+l\} & \mapsto M_{k,l}
\qquad \text{for $k,l \ge 1$},
\end{align}
permutations of the cells mapping to permutations of the generators
$M_n$ and $M_{k,l}$ of the $B_\infty$-operad. Finally define
\[
H_\bullet (J_3(n),J_2(n)) \to B_\infty(n), \qquad n \ge2,
\]
as zero.

Since $\dim \{1, \dots, n\} = n-2 = - \abs {M_n}$ and $\dim \{1, \dots,
k\}\{1, \dots, l\} = k+l-1 = - \abs {M_{k,l}}$ and the action of the
symmetric groups is respected by construction, the maps \eqref{corr1}
and \eqref{corr2} define a morphism $\oE \to B_\infty$ of graded
operads. The only thing which remains to be checked is the
compatibility of this morphism with the differentials. 

We will compute the boundary (differential) $d:=d^1$ on $\oE$. Each space
$J_i \setminus J_{i-1}$, $i=1,2,3$, is a disjoint union of
Fox-Neuwirth cells, which form a basis of the relative homology
$H_\bullet (J_i, J_{i-1})$. We will study the action of $d$ on this
basis.

\emph{$i=1$}.\ \
Let us start with $i=1$, when the points in a Fox-Neuwirth cell group
on a single vertical line. For $n \ge 2$ the boundary of the cell
$\{1,2, \dots , n\}$ in $\oE$ may be computed as follows:

\noindent
\let\picnaturalsize=N
\def\picsize{3in}
\def\picfilename{d32.eps}
\ifx\nopictures Y\else{\ifx\epsfloaded Y\else\input epsf \fi
\let\epsfloaded=Y
\centerline{\ifx\picnaturalsize N\epsfxsize \picsize\fi
\epsfbox{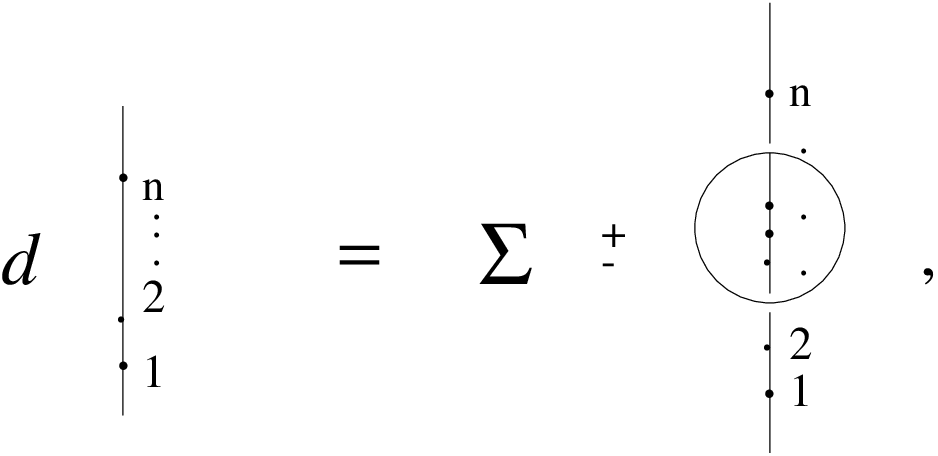}}}\fi

\noindent
where the left-hand side denotes the boundary of the cell $\{1, \dots,
n\}$ and the right-hand side denotes a linear combination of
Getzler-Jones cells obtained as operad compositions $\{1, \dots, k \}
\circ_i \{ 1, \dots, l\}$, $k, l \ge 2$, of two Fox-Neuwirth
cells. This equation turns into equation \eqref{A_infty}, where all
the terms with $i> 1$ or $j>1$ are moved to the right-hand side. The
signs here and henceforth in the proof are compatible with the signs
in \eqref{A_infty}--\eqref{Leibniz}, if the orientations on
Getzler-Jones cells are chosen as in the proof of
Proposition~\ref{prop1.2}.

\emph{$i=2$}.\ \ 
Cells for $i=2$ are configurations of points on two vertical
lines. The boundary of a cell $\{1, \dots k\} \{k+1, \dots, k+l\}$ may
be described as follows:
\smallskip

\noindent
\let\picnaturalsize=N
\def\picsize{4.6in}
\def\picfilename{d2n.eps}
\ifx\nopictures Y\else{\ifx\epsfloaded Y\else\input epsf \fi
\let\epsfloaded=Y
\centerline{\ifx\picnaturalsize N\epsfxsize \picsize\fi
\epsfbox{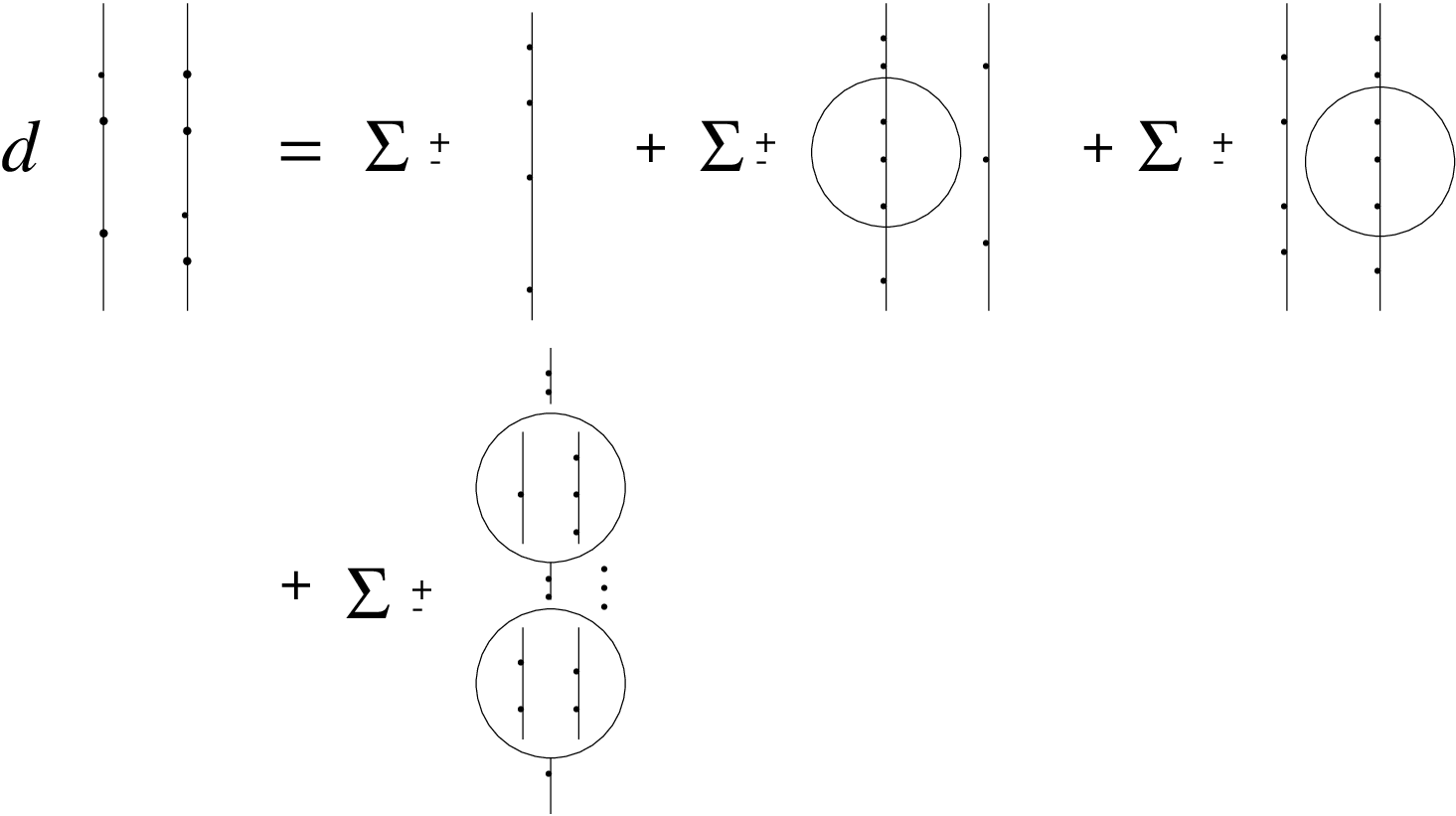}}}\fi

\noindent
This equation translates under the correspondence \eqref{corr1},
\eqref{corr2} into the identity \eqref{Leibniz} in which all terms but
those containing $d=M_1$ are moved to the right-hand side. Note that
the first sum on the figure corresponds to the terms $M_{k+l} (c_1,
\dots, c_{k+l})$, where $c_1, \dots, c_{k+l}$ is a shuffle of $\{a_1,
\dots, a_k\}$ and $\{b_1, \dots, b_l\}$, which show up on the
left-hand side of \eqref{Leibniz} for all pairs $(k_i,l_i)$ being
$(0,1)$ or $(1,0)$. The rest of the left-hand side of \eqref{Leibniz}
is the last term on the figure above.

\emph{$i=3$}.\ \ 
Cells for $i=3$ are configurations of points on at least three
different vertical lines. The boundary of a cell with at least four
vertical lines will produce the identity $0=0$ under the morphism $\oE
\to B_\infty$, because of a dimension argument: The boundary of such
cell has a dimension $ \ge n$, while (multiple) operad compositions of
Fox-Neuwirth cells from $J_2$ will have dimensions $\le
n-1$. Therefore, the boundary of a cell with at least four vertical
lines will have no terms which are compositions of Fox-Neuwirth cells
from $J_2$. Thus, the only nontrivial identity to be checked in
$B_\infty$ comes from the lowest Fox-Neuwirth cells in $J_3\setminus
J_2$, those made out of configurations of points on three vertical
lines.

The following figure describes the differential of such cell in $\oE$.

\let\picnaturalsize=N
\def\picsize{4.6in}
\def\picfilename{dmnk.eps}
\ifx\nopictures Y\else{\ifx\epsfloaded Y\else\input epsf \fi
\let\epsfloaded=Y
\centerline{\ifx\picnaturalsize N\epsfxsize \picsize\fi
\epsfbox{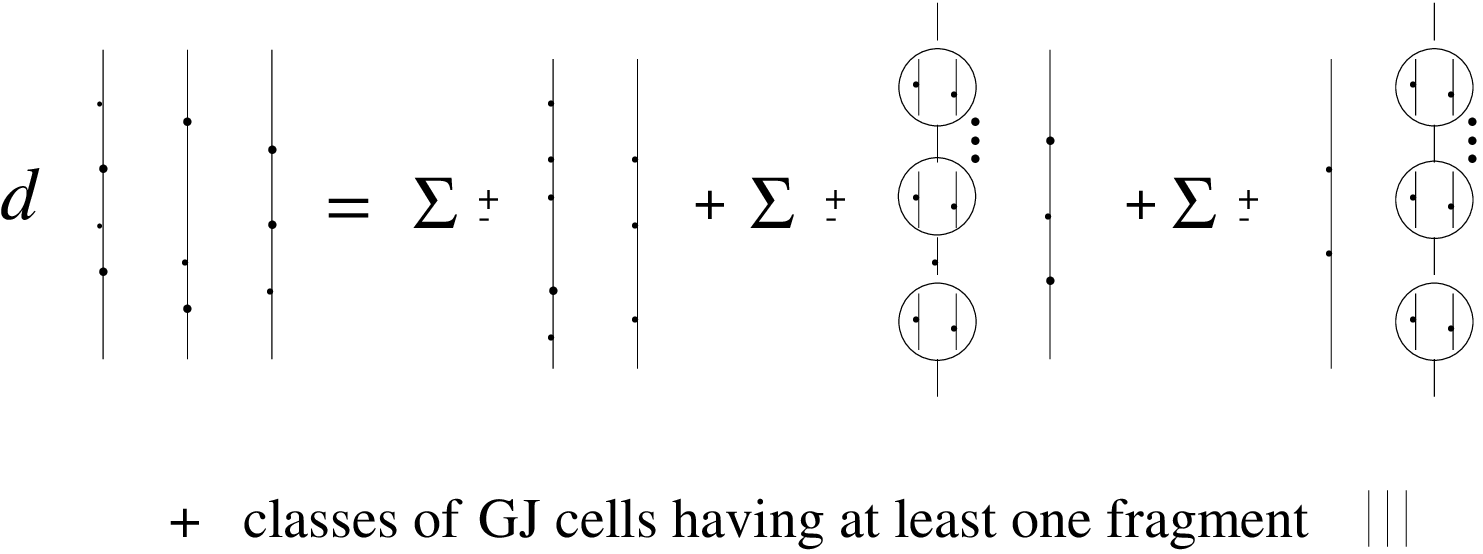}}}\fi
\smallskip

\noindent
Under the morphism $\oE \to B_\infty$, this identity turns into
\eqref{assoc} where all terms are moved to the right-hand side. Note
that the first sum on the figure corresponds to the terms in
\eqref{assoc} where all $(l_i,m_i)$ or all $(k_i,l_i)$ are either
$(0,1)$ or $(1,0)$.
\end{proof}

\begin{crl}
\label{formality}
The operad $\oE$ is formal, \emph{i.e.,} quasi-isomorphic to its
homology $e_2$. There is a morphism of operads $E^1 \to \oE$, unique
up to homotopy.
\end{crl}

\begin{proof}
In \cite{tamarkin:def} Tamarkin has constructed an operad morphism
$B_\infty \to e_2$, which is surjective on homology. Composing it with
our morphism $\oE \to B_\infty$, we get a morphism $\phi: \oE \to
e_2$. To prove the formality of $\oE$, it is enough to show that it is
a quasi-isomorphism. Moreover, it suffices to show that $\phi$ is
surjective on homology, because of a graded dimension argument.

Note that the operad $e_2$ is generated by $e_2(2)$, therefore we just
need to show that the second component $\phi(2): \oE(2) \to e_2(2)$ of
$\phi$ is surjective on homology. The morphism $\phi$ is a composition
$\oE \to B_\infty \to e_2$. According to Tamarkin \cite[Theorem
4.2.1]{tamarkin:def}, the induced homology morphism $H_\bullet
(B_\infty(2)) \to e_2(2)$ is an isomorphism. On the other hand, notice
that our morphism $\oE(2) \to B_\infty(2)$ is an isomorphism, because
$J_3 (2) = J_2 (2)$. In particular, it induces an isomorphism
$H_\bullet(\oE(2)) \to H_\bullet(B_\infty (2))$ of the homology. Thus,
the composition $H_\bullet (\phi): H_\bullet (\oE(2)) \to e_2(2)$ is
an isomorphism, which shows that $\oE$ is formal.

The existence of a unique up to homotopy morphism $E^1 \to \oE$
follows from the fact that both operads are quasi-isomorphic to $e_2$
and $E^1$ is a minimal model of $e_2$, see Corollary~\ref{crl}.
\end{proof}

\begin{quest}
Is the homology of the $B_\infty$-operad isomorphic to the G-operad
$e_2$? If yes, it will automatically be formal via Tamarkin's morphism
$B_\infty \to e_2$.
\end{quest}

\section{Action on the Hochschild complex}

\subsection{The Hochschild complex of an associative algebra}

Let us recall some notions related to the \h complex and some
properties of it. Let $A$ be an associative algebra and $C^n(A, A) =
\Hom ( A^{\otimes n} , A)$ its \emph{\h complex} with the \emph{\h
differential}:
\begin{multline}
\label{Hd}
(d x) (a_1, \dots, a_{n+1}) \\
\begin{aligned}
 := \; &
\; a_1 x(a_2, \dots, a_{n+1}) \\
 &  + \sum_{i=1}^n (-1)^i x (a_1, \dots , a_{i-1}, a_i a_{i+1},
a_{i+2}, \dots , a_{n+1})\\
 & -  (-1)^{n} x (a_1, \dots, a_n) a_{n+1},
\end{aligned}
\end{multline}
for $x \in C^n (A,A)$, $a_1, \dots, a_{n+1} \in A$.  The sign
$(-1)^{n}$ above is equal to $(-1)^{\abs x}$, $\abs x := n$ being the
\emph{degree} of the cochain $x$ in the \h complex $C^\bullet =
C^\bullet (A,A)$.

We will use the following operations on the Hochschild complex.  The
\emph{dot product} is defined as the usual cup product:
\begin{equation}
\label{dot-eq}
(x \cdot y) (a_1, \dots , a_{k+l}) = x(a_1, \dots, a_k) y (a_{k+1},
\dots , a_{k+l})
\end{equation}
for any $k$- and $l$-cochains $x$ and $y$ and $a_i \in A$. The
following collection of multilinear operations, called \emph{braces,}
cf.\ Kadeishvili \cite{kadeishvili} and Getzler \cite{g:cartan}, on
$C^\bullet$ is defined as
\begin{multline*}
\{x\} \{ x_1, \dots, x_n\} (a_1, \dots, a_m)
 := \\
  \sum
(-1)^\varepsilon x(a_1, \dots , a_{i_1},  x_1 (a_{i_1+1}, \dots), \dots,
a_{i_n},
  x_n(a_{i_n+1}, \dots),  \dots, a_m)
\end{multline*}
for $x , x_1, \dots , x_n \in C^\bullet$, $a_1, \dots, a_m \in A$,
where the summation runs over all possible substitutions of $x_1,
\dots, x_n$ into $x$ in the prescribed order and $\varepsilon :=
\sum_{p=1}^{n} (\abs {x_{p}} -1) i_{p} $. The braces $\{x\}\{ x_1,
\dots, x_n\}$ are homogeneous of degree $-n$, \emph{i.e}., $\abs
{\{x\}\{ x_1, \dots, x_n\}} = \abs x + \abs {x_1} + \dots + \abs {x_n}
- n$. We will also adopt the following notation:
\begin{align*}
x \circ y &:= \{x\} \{y\},\\
[x,y] & := x \circ y - (-1)^{(\abs {x}-1)(\abs {y}-1)} y \circ x.
\end{align*}
The \emph{G-bracket} $[x,y]$ defines the structure of a G-algebra on
the \h \coh\ $H^\bullet(A,A)$. The bracket was introduced by
Gerstenhaber \cite{gerst} in order to describe the obstruction for
extending a first-order deformation of the algebra $A$ to the second
order. The following definition of the bracket is due to Stasheff
\cite{jim:bracket}. Considering the tensor coalgebra $T(A) =
\bigoplus_{n=0}^\infty A^{\otimes n}$ with the comultiplication
$\Delta (a_1 \otimes \dots \otimes a_n) = \sum_{k=0}^n (a_1 \otimes
\dots \otimes a_k) \otimes (a_{k+1} \otimes \dots \otimes a_n)$, we
can identify the \hc s $\Hom (A^{\otimes n}, A)$ with the graded
coderivations $\coder T(A)$ of the tensor coalgebra $T(A)$. Then the
G-bracket $[x,y]$ is defined as the (graded) commutator of
coderivations. In fact, the \h complex $C^\bullet$ is a differential
graded Lie algebra with respect to this bracket.

In addition, the dot product \eqref{dot-eq} and the Hochschild
differential \eqref{Hd} define the structure of a DG associative
algebra on $C^\bullet (A,A)$.

\subsection{The structure of a $B_\infty$-algebra on the Hochschild complex}

Define the structure of a $B_\infty$-algebra on $C^\bullet (A,A)$ as
follows:
\begin{align*}
M_0 & : = 0,\\
M_1 & := d,\\
M_2 (x_1,x_2) & := x_1 \cdot x_2,\\
M_n & := 0 \qquad \text{for $n > 2$},\\
M_{0,1} = M_{1,0} & := \id,\\
M_{0,n} = M_{n,0} & := 0 \qquad \text{for $n > 1$},\\
M_{1,n} (x; x_1, \dots, x_n) & := \{x\}\{x_1, \dots, x_n\}
\qquad \text{for $n \ge 0$},\\
M_{k,l} & := 0 \qquad \text{for $k > 1$},
\end{align*}
where $x, x_1, \dots, x_n \in C^\bullet (A,A)$.

\begin{thm}
\label{B_infty}
These operations define the structure of a $B_\infty$-algebra on the
\h complex $C^\bullet$.
\end{thm}

\begin{rem}
The braces were defined by Kadeishvili \cite{kadeishvili} and Getzler
\cite{g:cartan}, the identities among them and the dot product were
written down in \cite{gv1}, where this structure was called a homotopy
G-algebra structure. The fact that this algebraic data defines a
$B_\infty$-structure was noticed by Getzler and Jones in \cite{gj}.
\end{rem}

\begin{proof}

Taking into account the vanishing operations $M_n$ and $M_{k,l}$ and
rewriting the rest in terms of the dot product and braces, the
identities \eqref{A_infty} through \eqref{Leibniz} can be simplified
as follows.

The identities \eqref{A_infty} for $n=1$, 2, and 3, are equivalent to
\begin{gather}
\label{onthenose}
d^2 =  0, \\
d (x_1 x_2) = (dx_1)x_2 + (-1)^{\abs {x_1}} x_1 dx_2,\\
(x_1x_2) x_3 = x_1 (x_2 x_3),
\end{gather}
respectively.

The identities \eqref{assoc} are nontrivial only for $k=1$, when they
turn into the following.
\begin{multline}
\label{MM}
\sum_{0 \le i_1 \le \dots \le i_l \le m} (-1)^{\epsilon} \{x\} \{ z_1, \dots,
z_{i_1}, \{y_1\} \{z_{i_1+1}, \dots \}, \dots , \\
z_{i_l}, \{ y_{l}\}
\{y_{i_l + 1}, \dots\}, \dots , y_m\}
\\
= \{ \{ x\} \{y_1, \dots, y_l\}\} \{ z_1, \dots, z_m\},
\end{multline}
where $\epsilon = \sum_{p=1}^l (\abs{y_p}-1)
\sum_{q=1}^{i_p}(\abs{z_q} -1)$.

The identities \eqref{Leibniz} are nontrivial only when $k=1$ and 2. For
$k=1$ they rewrite as the following family of identities:
\begin{multline}
\label{k=1}
\begin{aligned}
d \{x\} \{y_1, \dots, y_{l}\} - (-1)^{\abs{x} (\abs{y_1}-1)} y_1 \cdot
 \{ x \} \{y_2, \dots, y_l \} \\ + (-1)^{\abs{x} + \abs{y_1} + \dots +
 \abs{y_{l-1}} + l -1} \{x \} \{y_1 , \dots, y_{l-1}\} \cdot y_l
\end{aligned}
\\
= \{ dx \} \{ y_1 , \dots, y_l\} \\
- \sum_{i=0}^{l-1} (-1)^{|x| + \abs{y_1} + \dots + \abs{y_i} - i} \{ x \} \{ y_1 , \dots, y_i, d
y_{i+1}, \dots, y_l \}
\\
- \sum_{i=0}^{l-2} (-1)^{|x| + \abs{y_1} + \dots + \abs{y_{i+1}} - i}
 \{ x \} \{ y_1 , \dots, y_i, y_{i+1} \cdot y_{i+2}, \dots, y_l\},
\end{multline}
for each $l \ge 1$. For $k=2$, Equations~\eqref{Leibniz} turn into
\begin{multline}
\label{k=2}
\sum_{0 \le l_1 \le l} (-1)^{\abs{x_2}(\abs{y_1} + \dots
+\abs{y_{l_1}} - l_1)} (\{x_1\} \{y_1, \dots, y_{l_1}\}) \cdot (\{ x_2
\} \{y_{l_1+1}, \dots, y_l \} )
\\
= \{ x_1 \cdot x_2 \} \{ y_1 , \dots, y_l\},
\end{multline}
for each $l \ge 1$.

All these identities for the operations on the \h complex may be
checked directly. Some of the identities are classical, see
\emph{e.g}., Gerstenhaber \cite{gerst}, the others were not noticed
until more recently, see \cite{gv1}. One can find a detailed
verification of the identities in Khalkhali's paper \cite{kh}.
\end{proof}

Combining Theorems \ref{thm:relation} and \ref{B_infty} and
Corollary~\ref{formality}, we come to the following solution of
Deligne's Conjecture.
\begin{crl}
\label{final}
The operad morphism $\oE \to B_\infty$ and the above action of
$B_\infty$ on the \h complex $C^\bullet$ define on $C^\bullet$ the
natural structure of an algebra over the operad $\oE$, which is
quasi-isomorphic to its homology $e_2$.
\end{crl}

\begin{rem}
This corollary along with Corollary~\ref{formality} also yields the
natural structure of a $G_\infty$-algebra on $C^\bullet$, recovering a
result of Tamarkin \cite{tamarkin:def}.
\end{rem}

\begin{rem}
A complex of vector spaces with operations $x_1 \cdot x_2$ and
$\{x\}\{x_1, \dots, x_n\}$ for $n \ge 0$ satisfying identities
\eqref{onthenose}--\eqref{k=2} was called a \emph{homotopy G-algebra}
in \cite{gv1,gv2}. Kadeishvili rediscovered the same notion in
\cite{kadeishvili:hirsch} under the name of an \emph{associative
Hirsch algebra}.
\end{rem}

\section{Correction of \cite{kvz}}

As we have already mentioned, the paper \cite{kvz} of Kimura,
Zuckerman, and the author used the notion of a $G_\infty$-algebra as
an algebra over the Getzler-Jones' cellular operad, which was later
noticed not to be cellular. The following changes have to be made to
correct the resulting error.

In Section 4, after describing Getzler-Jones' cellular partition
$K_\bullet \FF$, we should emphasize that it is not a cellular
operad. However, there is a way to produce a DG operad out of
it. Namely, notice that as a graded operad, $K_\bullet \FF$ is free on
the collection of Fox-Neuwirth cells. The problem is that the
differential is not well defined, in general. For those Fox-Neuwirth
cells $C_1$ whose geometric boundary is the union of Getzler-Jones
cells, the differential $dC_1$ is defined as the boundary operator. For
those Fox-Neuwirth cells $C_2$ whose geometric boundary is not the
union of Getzler-Jones cells, define the differential as a new
generator $dC_2$. Take the free graded operad $K_\bullet$ generated by
the Fox-Neuwirth cells and the differentials of the Fox-Neuwirth cells
of the second type. The differential now is well defined on the
generators: For a Fox-Neuwirth cell $C_1$ of the first type, the
differential $dC_1$ is a linear combination of Getzler-Jones cells;
for a Fox-Neuwirth cell $C_2$ of the second type, the differential is
the generator $dC_2$; finally $d(dC_2) = 0$. This differential extends
uniquely to the free graded operad $K_\bullet$. Definition~4.1 in
\cite{kvz} must be replaced by the following one.

\begin{df}
The \emph{weak $G_\infty$-operad} is the DG operad $K_\bullet$
constructed above. An algebra over it is called a \emph{weak
$G_\infty$-algebra}.
\end{df}

Every occurrence of the word ``$G_\infty$-algebra'' in \cite{kvz} must
be replaced with the words ``weak $G_\infty$-algebra''. Whenever the
``operad'' $K_\bullet \FF$ occurs in the sequel therein, it must be
replaced with the above operad $K_\bullet$. All the lower identities
written down in Sections 4.1 and 4.2 of \cite{kvz} are satisfied in a
weak $G_\infty$-algebra and do not require corrections, because the
Fox-Neuwirth cells used there are of the first type. For the same
reason, the $A_\infty$- and $L_\infty$-operads map naturally to the
weak $G_\infty$-operad, therefore, a weak $G_\infty$-algebra is
naturally an $A_\infty$- and $L_\infty$-algebra. With the replacement
of the $G_\infty$-algebras by weak $G_\infty$-algebras, all the
results of the paper are correct with the same proofs. Moreover,
Conjecture~2.3 of \cite{kvz} with the above change has been proven in
the meantime by Yi-Zhi Huang and Wenhua Zhao \cite{huang-zhao}. Thus
it must be renamed to a theorem, as follows.

\begin{thm}[Huang and Zhao]
Let $V^\bullet [\cdot]$ be a TVOA satisfying $G(0)^2=0$. Then the dot
product and the skew-symmetrization of the bracket defined in
\cite{kvz} can be extended to the structure of a weak
$G_\infty$-algebra on a certain completion of $V^\bullet [\cdot]$.
\end{thm}

\bibliographystyle{amsalpha}
\bibliography{../dennis/op}
\end{document}